\documentclass[12pt]{article}
\usepackage{graphicx, color}
\usepackage{url}
\usepackage{amsmath, amsfonts, rotating}
\usepackage{amsmath}
\usepackage{amsfonts}
\usepackage{amssymb}
\usepackage[ot2enc]{inputenc}
\usepackage[russian,english]{babel}

\title{Iterated partial summations applied to finite-support discrete distributions}
\author{Michaela Ko\v{s}\v{c}ov\'a, Radoslav Harman, J\'an Ma\v{c}utek\\ \\Department of Applied Mathematics and Statistics,\\ Comenius University in Bratislava, Slovakia\\ \\ e-mails: michaela.koscova@fmph.uniba.sk, \\ radoslav.harman@fmph.uniba.sk, jmacutek@yahoo.com}

\begin{document}

\maketitle
\abstract{The problem of iterated partial summations is solved for some discrete distributions defined on discrete supports. The power method, usually used as a computational approach to finding matrix eigenvalues and eigenvectors, is in some cases an effective tool to prove the existence of the limit distribution, which is then expressed as a solution of a system of linear equations. Some examples are presented. \\}


\section{Introduction}\label{sec:matrix_notation}

Let $\{P_j^*\}_{j=0}^\infty$ and $\{P_j\}_{j=0}^\infty$ be discrete probability distributions defined on the set of non-negative integers. A general form of partial-sums distributions was introduced in \cite{macutek_2003}. The distribution $\{P_j\}_{j=0}^\infty$ is the result of a partial summation applied to $\{P_j^*\}_{j=0}^\infty$ if
\begin{equation}
P_x=c\sum_{j=x}^\infty g(j)P_j^*,\quad\quad x=0,1,2,\dots,\label{eq:partial_summation}
\end{equation}
where $g(j)$ is a real function and $c$ a normalization constant (which ensures that the sequence $\{P_j\}_{j=0}^\infty$ is a proper probability distribution, i.e., it sums to 1). The distributions $\{P_j^*\}_{j=0}^\infty$ and $\{P_j\}_{j=0}^\infty$ are called parent and descendant, respectively. Some special cases of (\ref{eq:partial_summation}) are mentioned also in the comprehensive monograph on discrete distributions \cite{johnson_et_al_2005}, pp.~508-512.

Hereafter, we restrict our considerations to parent distributions $\{P_0^*,P_1^*,P_2^*,\dots, P_{S-1}^*\}$, i.e. to discrete distributions defined on a finite support of the size $S$. As the probabilities $P^*_j$ in (\ref{eq:partial_summation}) are zero for $j\geq S$ in this case, the partial summation (\ref{eq:partial_summation}) can be written as
\begin{equation*}
P_x=c\sum_{j=x}^{S-1}g(j)P_j^*,\quad\quad x=0,1,\dots,S-1 ,
\end{equation*}
or equivalently as
\begin{equation}
\mathbb{P}=cA\mathbb{P^*},\label{eq:matrix_equation}
\end{equation}
where $\mathbb{P}$ and $\mathbb{P^*}$ are the vectors of probabilities $\{P_0,P_1,\dots,P_{S-1}\}^\top$ and $\{P_0^*,P_1^*,\dots,P_{S-1}^*\}^\top$, respectively. Matrix $A$ is of dimension $S\times S$, with the following structure:
\begin{equation}
A=\left(
\begin{matrix}
    g(0) & g(1) & g(2) & g(3) &\dots & g(S-1)  \\
    0 & g(1) & g(2) & g(3) & \dots & g(S-1)  \\
    0 & 0 & g(2) & g(3) & \dots & g(S-1)  \\
    0 & 0 & 0 & g(3) & \dots & g(S-1)  \\
    \vdots & \vdots & \vdots & \vdots & \ddots & \vdots\\
    0 & 0 & 0 & 0 & \dots & g(S-1)
\end{matrix}
\right).\label{eq:matrix_A}
\end{equation}

\section{Iterated partial summations}\label{sec:iterated}

The partial summation (\ref{eq:partial_summation}) can be applied iteratively. The descendant distribution becomes a parent of another distribution, i.e.
\begin{alignat*}{2}
P_x^{(1)}=&c_1\sum_{j=x}^{S-1} g(j)P_j^*,\quad\quad x=0,1,\dots,{S-1} \ ,
\nonumber\\
P_x^{(2)}=&c_2\sum_{j=x}^{S-1} g(j)P_j^{(1)},\quad\quad x=0,1,\dots,{S-1} \ ,\nonumber\\
\quad\quad\vdots\nonumber\\
P_x^{(n)}=&c_n\sum_{j=x}^{S-1} g(j)P_j^{(n-1)},\quad\quad x=0,1,\dots,{S-1} \ ,\nonumber\\
\quad\quad\vdots\nonumber
\end{alignat*}
with $c_i$, $i=1,2,3,\dots \ $ being normalization constants. The distribution $\{P_x^*\}_{x=0}^{S-1}$ will be called the original parent. We will now investigate properties of the sequence of the descendant distributions, especially the question under which conditions the limit of this sequence exists. The existence of the limit for iterated partial summations applied to discrete distributions with infinite supports for a constant function $g(j)$ was proved in \cite{macutek_2006}.
%begin{equation}
%\lim_{i\rightarrow\infty}P_x^{(i)},\quad\quad x=0,1,\dots,S-1.
%\end{equation}

In the following, we will not consider the normalization constants. Then the matrix notation (see (\ref{eq:matrix_equation})) of the iterated partial summations is
\begin{alignat*}{2}
\mathbb{Q}^{(1)}=&A\mathbb{P^*},\nonumber\\
\mathbb{Q}^{(2)}=&A\mathbb{Q}^{(1)}=AA\mathbb{P^*}=A^2\mathbb{P^*},\nonumber\\
\mathbb{Q}^{(3)}=&A\mathbb{Q}^{(2)}=AA^2\mathbb{P^*}=A^3\mathbb{P^*},\nonumber\\
\quad\quad\vdots\nonumber\\
\mathbb{Q}^{(n)}=&A\mathbb{Q}^{(n-1)}=A^n\mathbb{P^*},\nonumber\\
\quad\quad\vdots\nonumber
\end{alignat*}
The $i$-th descendant probability distribution can be obtained by the normalization of the vector $\mathbb{Q}^{(i)}=(Q_0^{(i)},Q_1^{(i)},\dots,Q_{S-1}^{(i)})^\top$. 

Denote $\Vert u \Vert_1$, $\Vert u \Vert_2$ the L1-norm and the L2-norm of vector $u$, respectively. If the limit of the sequence of the descendant distributions exists, it can be written as
\begin{equation}
\mathbb{P^{(\infty)}}=\lim_{n\rightarrow\infty}\frac{\mathbb{Q}^{(n)}}{\Vert\mathbb{Q}^{(n)}\Vert_1}=\lim_{n\rightarrow\infty}\frac{A^n\mathbb{P^*}}{\Vert A^n\mathbb{P^*}\Vert_1}.\label{eq:limit}
\end{equation}

\section{Application of the power method}\label{sec:power_method}

The power method is one of computational approaches to finding matrix eigenvalues (see e.g. \cite{golub_vanloan_1996}, pp.~330-332). We apply it to matrix $A$ (denote its eigenvalues by $\lambda_0, \dots, \lambda_{S-1}$; we remind that in general they need not be distinct) and vector $\mathbb{P^*}$ from (\ref{eq:matrix_equation}). To satisfy the conditions of the method, suppose that $A$ is diagonalizable and that it has a unique dominant eigenvalue $\lambda_k$ (i.e., there exists $k$ such that $\vert\lambda_k \vert > \vert\lambda_i \vert$, $i\neq k$).

If all $S$ elements of $\mathbb{P^*}$ are non-zero and if $\mathbb{P^*}$ is not a non-dominant eigenvector of $A$, then
\begin{equation*}
\lim_{n\rightarrow\infty}\frac{A^n\mathbb{P^*}}{\Vert A^n\mathbb{P^*}\Vert_2}=v,
\end{equation*}
where $v$ is the dominant eigenvector of $A$ (i.e., the one which corresponds to the dominant eigenvalue). Under these conditions, the power methods implies the existence of $\lim_{n\rightarrow\infty}\mathbb{P}^{(n)}$, see (\ref{eq:limit}), with
\begin{equation*}
\mathbb{P^{(\infty)}}=\lim_{n\rightarrow\infty}\mathbb{P}^{(n)}=\frac{v}{\Vert v\Vert_1}.
\end{equation*}
The matrix $A$ from (\ref{eq:matrix_A}) is an upper triangular matrix, which means that its eigenvalues are its diagonal entries, i.e.
\begin{equation*}
\lambda_j=g(j),\quad\quad j=0,1,\dots,S-1.
\end{equation*} 
Consequently, to determine the dominant eigenvalue $\lambda_D$ of $A$ it is necessary to find
\begin{equation*}
D=
%\arg\max_{j\in\{0,1,\dots,S-1\}}\vert\lambda_j\vert
\arg\max_{j\in\{0,1,\dots,S-1\}}\vert g(j)\vert. \label{dominant_eigenvalue}
\end{equation*}

Let $D=k$, i.e., let the dominant eigenvalue be $\lambda_k=g(k)$. %Additionally, let the dominant eigenvalue be unique. 
The eigenvector corresponding to the dominant eigenvalue $\lambda_k$ is the solution of the of linear equations
\begin{equation*}
Av=\lambda_kv,
\end{equation*}
or, equivalently,
\begin{equation*}
(A-\lambda_kI)v=0.
\end{equation*}
For matrix $A$ from (\ref{eq:matrix_equation}) we obtain
\newpage
\begin{sideways}\begin{minipage}{\textheight} 		
\begin{equation*}
\left(
\begin{matrix}
    g(0)-g(k) & g(1) & g(2)  & \dots & g(k-1) & g(k) & g(k+1)  &\dots & g(S-2) & g(S-1)  \\
    0  & g(1)-g(k) & g(2)  & \dots & g(k-1) & g(k) & g(k+1)  &\dots & g(S-2) & g(S-1)  \\
    \vdots & \vdots  &  \vdots &  \ddots & \vdots & \vdots & \vdots &  \ddots & \vdots & \vdots\\
    0 & 0 & 0  & \dots & g(k-1)-g(k) & g(k) & g(k+1) & \dots & g(S-2) & g(S-1)\\
    0 & 0 & 0  & \dots & 0 & 0 & g(k+1) & \dots & g(S-2) & g(S-1)\\
    0 & 0 & 0  & \dots & 0 & 0 & g(k+1)-g(k) & \dots & g(S-2) & g(S-1)\\
    \vdots   & \vdots & \vdots &  \ddots & \vdots & \vdots & \vdots &  \ddots & \vdots & \vdots\\
    0 & 0 & 0  & \dots & 0 & 0 & 0 & \dots & g(S-2)-g(k) & g(S-1)\\
    0 & 0 & 0  & \dots & 0 & 0 & 0 & \dots & 0 & g(S-1)-g(k)
\end{matrix}
\right)
\left(\begin{matrix}
    v_0    \\
    v_1    \\
    \vdots    \\
    v_{k-1}    \\
    v_k\\
    v_{k+1}\\
    \vdots \\
    v_{S-2}\\
    v_{S-1}
\end{matrix}\right)
=0.
\end{equation*}
\end{minipage}\end{sideways}
\newpage
This system of linear equations yields the solution
\begin{equation}\label{eq:eigenvector}
v=
%\left(\begin{matrix}
%    v_0    \\
%    v_1    \\
%    v_2\\
%    \vdots    \\
%    v_k\\
%    v_{k+1}\\
%    \vdots \\
%    v_{S-1}
%\end{matrix}\right)=
\left(\begin{matrix}
    t    \\
    t\left(1-\frac{g(0)}{g(k)}\right)    \\
    t\left(1-\frac{g(0)}{g(k)}\right)\left(1-\frac{g(1)}{g(k)}\right)    \\
    \vdots    \\
    t\prod_{j=1}^k\left(1-\frac{g(j-1)}{g(k)}\right)\\
    0\\
    \vdots \\
    0
\end{matrix}\right), \quad t\in\mathbb{R}.
\end{equation}

\section{Example: Iterated Katz partial summations}\label{sec:example}

Discrete distribution $\{R_j^*\}_{j=0}^{n/\infty}$ belongs to the Katz family (see e.g. \cite{wimmer_altmann_1999}, pp.~324-325) with the parameters $\alpha\geq0$, $\beta<1$ if  
\begin{equation*}
\frac{R_{x+1}}{R_x}=\frac{\alpha+\beta x}{x+1},\quad\quad x=0,1,2,\dots~.
\end{equation*}
The Katz partial summation, i.e. the summation 
\begin{equation*}
P_x=\sum_{j=x}^\infty g(j)P_j^*,\quad\quad x=0,1,2,\dots\label{eq:partial_summation1}
\end{equation*}
with
\begin{equation}
g(j)=1-\frac{R_{j+1}}{R_j}=1-\frac{\alpha+\beta j}{j+1}=\frac{(1-\alpha)+(1-\beta)j}{j+1},\quad\quad j=0,1,2,\dots\label{eq:g_Katz}
\end{equation}
was analyzed in \cite{macutek_2001}.

Consider a finite-support discrete distribution $\{P_0^*,P_1^*,P_2^*,\dots, P_{S-1}^*\}$. Function $g(j)$ from (\ref{eq:g_Katz}) is increasing if $\alpha > \beta$, constant if $\alpha = \beta$, and decreasing if $\alpha < \beta$. Thus, if $\alpha \neq \beta$, all eigenvalues of matrix $A$ are distinct, which is a sufficient condition for its diagonalizability (see \cite{horn_johnson_2013}). If, in addition, $\vert g(0)\vert \neq \vert g(S-1)\vert$, there exists the unique dominant eigenvalue and the power method can be applied. The dominant eigenvector $v$ from (\ref{eq:eigenvector}) can be expressed as
\begin{equation*}
v=\left(\begin{matrix}
    v_0    \\
    v_1    \\
    v_2\\
    \vdots    \\
    v_k\\
    v_{k+1}\\
    \vdots \\
    v_{S-1}
\end{matrix}\right)=
\left(\begin{matrix}
    t    \\
    t\frac{\alpha-\beta}{(1-\alpha)+(1-\beta)k}k    \\
    t\left(\frac{\alpha-\beta}{(1-\alpha)+(1-\beta)k}\right)^2\frac{k(k-1)}{2}    \\
    \vdots    \\
     t\left(\frac{\alpha-\beta}{(1-\alpha)+(1-\beta)k}\right)^k\frac{k!}{k!}\\
    0\\
    \vdots \\
    0
\end{matrix}\right), \quad t\in\mathbb{R}.
\end{equation*}
We use the parameter $t$ to scale the vector $v$ so that the sum of the vector elements is equal to 1, i.e.
\begin{equation*}
t=\left(\frac{(1-\alpha)+(1-\beta)k}{(1-\beta)(k+1)}\right)^k.
\end{equation*}
Therefore,
\begin{equation*}
\mathbb{P}^{(\infty)}=
\left(\begin{matrix}
    {P}^{(\infty)}_{0}    \\
    {P}^{(\infty)}_{1}    \\
    {P}^{(\infty)}_{2}	\\
    \vdots    \\
    {P}^{(\infty)}_{k}\\
    {P}^{(\infty)}_{k+1}\\
    \vdots \\
    {P}^{(\infty)}_{S-1}
\end{matrix}\right)=
\left(\begin{matrix}
    \left(1-\frac{\alpha-\beta}{(1-\beta)(k+1)}\right)^k   \\
    k\frac{\alpha-\beta}{(1-\beta)(k+1)}\left(1-\frac{\alpha-\beta}{(1-\beta)(k+1)}\right)^{k-1}    \\
    {k\choose 2}\left(\frac{\alpha-\beta}{(1-\beta)(k+1)}\right)^2\left(1-\frac{\alpha-\beta}{(1-\beta)(k+1)}\right)^{k-2}      \\
    \vdots    \\
     \left(\frac{\alpha-\beta}{(1-\beta)(k+1)}\right)^k\\
    0\\
    \vdots \\
    0
\end{matrix}\right),
\end{equation*}
which means $\mathbb{P}^{(\infty)}\sim Bin\left(k;\frac{\alpha-\beta}{(1-\beta)(k+1)}\right).$

However, as the function $g(j)$ from (\ref{eq:g_Katz}) is strictly monotonic in $j$ if $\alpha \neq \beta$, there are only two possible values of $k$, either 0 or $S-1$. If $k=S-1$, the iterated partial summations have the limit which is the binomial distribution with parameters $S-1$ and $\frac{\alpha-\beta}{(1-\beta)S}$. On the other hand, if $k=0$, the distribution $\mathbb{P}^{(\infty)}$ degenerates to the deterministic distribution. 

To inspect whether $k$ is equal to $0$ or to $S-1$ for a particular choice of parameters of the Katz family - i.e. for the parameters which appear in function $g(j)$ from (\ref{eq:g_Katz}) - it is sufficient to compare the values of $\vert g(0)\vert$ and $\vert g(S-1)\vert$, which means to solve the inequalities
%\begin{equation}
%\vert g(0)\vert\lesseqqgtr\vert g(S-1)\vert,
%\end{equation}
%which is the same as
\begin{equation}
\vert 1-\alpha\vert\lesseqqgtr\left\vert\frac{(1-\alpha)+(1-\beta)(S-1)}{S}\right\vert.\label{inequality}
\end{equation}
The result is shown in Figure \ref{fig:diagram}.

\begin{figure}[h]
\caption{The solution of the inequality (\ref{inequality})}
\begin{equation}
\begin{array}{llllllllll}
&\alpha=\beta&\Rightarrow	&	\text{?}	&&	&&&&\\
\nearrow&&&&&&&&&\\
\rightarrow&\alpha<\beta&\Rightarrow	&	\text{Det}	&&&	&&&\\
\searrow&&&&&&&&&\\
&\alpha>\beta&\rightarrow	&\alpha\leq 1	&	\Rightarrow&\text{Bin}&&&&\\
&&\searrow	&&	&&&&&\\
&& &\alpha>1	&\rightarrow	&\alpha+\beta-2\geq0&\Rightarrow&\text{Det}&&\\
&& &	&&&&&&\\
&& &	&	&&&S=\frac{\alpha-\beta}{2-\alpha-\beta}&\Rightarrow&\text{?}\\
&& &	&\searrow		&&\nearrow&&&\\
&& &	&	&\alpha+\beta-2<0&\rightarrow&S<\frac{\alpha-\beta}{2-\alpha-\beta}&\Rightarrow&\text{Det}\\
&& &	&	&&\searrow&&&\\
&& &	&	&&&S>\frac{\alpha-\beta}{2-\alpha-\beta}&\Rightarrow&\text{Bin}\\
\end{array}\nonumber
\end{equation}\label{fig:diagram}
\end{figure}

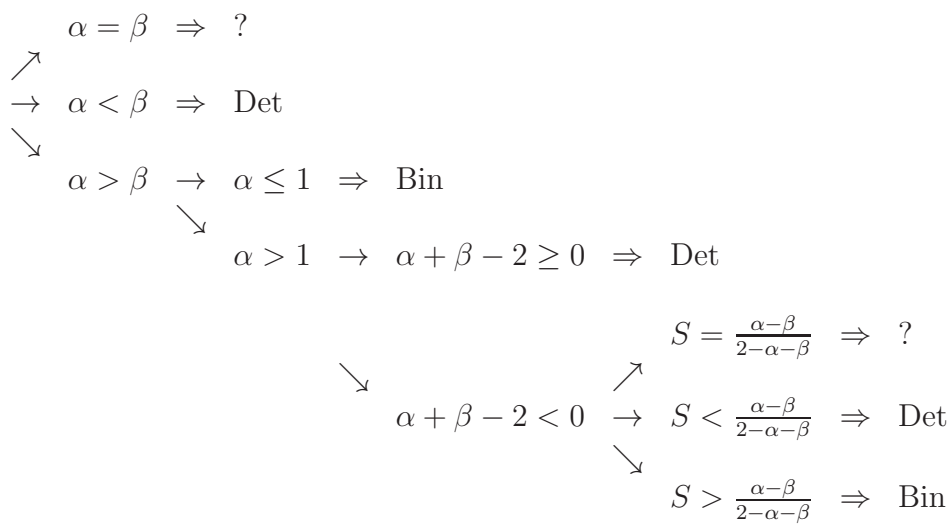
Figure \ref{fig:Katz_par_space} depicts the parametric space of the Katz partial summations. The iterated Katz partial summations with parameters from the green area result in the deterministic distribution. Those with parameters from the blue area result in the binomial distribution $Bin\left(S-1;\frac{\alpha-\beta}{S(1-\beta)}\right)$. We remind that these results are valid regardless of the original parent distribution $\{P_0^*,P_1^*,P_2^*,\dots, P_{S-1}^*\}$. The limit behaviour of the iterated Katz partial summations with parameters from the yellow area depends also on the size of the support of the original parent. If $S=\frac{\alpha-\beta}{2-\alpha-\beta}$, see Figure \ref{fig:diagram}, the limit distribution remains an open question (as the power method cannot be applied in this case). 

\setlength{\unitlength}{1cm}
\begin{figure}
\caption{The parametric space of Katz family}
\label{fig:Katz_par_space}
\begin{picture}(10,12.5)
%\color{brown}
\put(0.5,0){\vector(0,1){11.5}}
\put(0,5.5){\vector(1,0){13.5}}
\color{green}
\put(5,10.5){\line(1,-1){0.25}}	%nenakreslilo
\put(4.5,10.5){\line(1,-1){0.5}}
\put(4,10.5){\line(1,-1){0.75}}
\put(3.5,10.5){\line(1,-1){1}}
\put(3,10.5){\line(1,-1){1.25}}
\put(2.5,10.5){\line(1,-1){1.5}}
\put(2,10.5){\line(1,-1){1.75}}
\put(1.5,10.5){\line(1,-1){2}}
\put(1,10.5){\line(1,-1){2.25}}
\put(0.5,10.5){\line(1,-1){2.5}}
\put(0.5,10){\line(1,-1){2.25}}
\put(0.5,9.5){\line(1,-1){2}}
\put(0.5,9){\line(1,-1){1.75}}
\put(0.5,8.5){\line(1,-1){1.5}}
\put(0.5,8){\line(1,-1){1.25}}
\put(0.5,7.5){\line(1,-1){1}}
\put(0.5,7){\line(1,-1){0.75}}
\put(0.5,6.5){\line(1,-1){0.5}}
\put(0.5,6){\line(1,-1){0.25}}	%nenakreslilo
\color{red}
%\linethickness{1mm}
\thicklines
\put(0.5,5.5){\line(1,1){5}}
\thinlines
%\multiput(0.5, 5.5)(0.6, 0.6 ){11}{\line(1,1){0.4}}
\color{blue}
\multiput(0.5, 5)(0, -0.5 ){11}{\line(1,1){5}}
\put(1,0){\line(1,1){4.5}}
\put(1.5,0){\line(1,1){4}}
\put(2,0){\line(1,1){3.5}}
\put(2.5,0){\line(1,1){3}}
\put(3,0){\line(1,1){2.5}}
\put(3.5,0){\line(1,1){2}}
\put(4,0){\line(1,1){1.5}}
\put(4.5,0){\line(1,1){1}}
\put(5,0){\line(1,1){0.5}}
%\thinlines
\color{black}
\multiput(0.5, 10.5)(0.4, 0 ){33}{\line(1,0){0.2}}
%\put(0.5,10.5){\line(1,0){13}}
\put(5.5,0){\line(0,1){10.5}}
\put(5.5,10.5){\line(1,-1){8}}
\color{green}
\put(6,10.5){\line(1,-1){7.5}}
\put(6.5,10.5){\line(1,-1){7}}
\put(7,10.5){\line(1,-1){6.5}}
\put(7.5,10.5){\line(1,-1){6}}
\put(8,10.5){\line(1,-1){5.5}}
\put(8.5,10.5){\line(1,-1){5}}
\put(9,10.5){\line(1,-1){4.5}}
\put(9.5,10.5){\line(1,-1){4}}
\put(10,10.5){\line(1,-1){3.5}}
\put(10.5,10.5){\line(1,-1){3}}
\put(11,10.5){\line(1,-1){2.5}}
\put(11.5,10.5){\line(1,-1){2}}
\put(12,10.5){\line(1,-1){1.5}}
\put(12.5,10.5){\line(1,-1){1}}
\put(13,10.5){\line(1,-1){0.5}}
\color{yellow}
\put(6,0){\line(0,1){10}}
\put(6.5,0){\line(0,1){9.5}}
\put(7,0){\line(0,1){9}}
\put(7.5,0){\line(0,1){8.5}}
\put(8,0){\line(0,1){8}}
\put(8.5,0){\line(0,1){7.5}}
\put(9,0){\line(0,1){7}}
\put(9.5,0){\line(0,1){6.5}}
\put(10,0){\line(0,1){6}}
\put(10.5,0){\line(0,1){5.5}}
\put(11,0){\line(0,1){5}}
\put(11.5,0){\line(0,1){4.5}}
\put(12,0){\line(0,1){4}}
\put(12.5,0){\line(0,1){3.5}}
\put(13,0){\line(0,1){3}}
\color{black}
\put(13.5,5.2){\scriptsize{$\alpha$}}
\put(0.3,11.5){\scriptsize{$\beta$}}
%\put(0.5,5.5){\circle*{0.1}} %prazdny
\put(0.3,5.2){\scriptsize0}
\put(0.3,10.4){\scriptsize1}
%\put(5.5,5.5){\circle*{0.1}} 
\put(5.3,5.2){\scriptsize1}
\put(10.3,5.2){\scriptsize2}
%\put(0.5,10.5){\circle*{0.1}} 
%\put(0.6,10.6){\scriptsize{$A$}}
%\put(5.5,10.5){\circle*{0.1}} 		%prazdny treba
%\put(5.6,10.3){\scriptsize{$P$}}
%\put(0.5,0.5){\circle*{0.1}} 
%\put(0.6,0.3){\scriptsize{$G$}}
%\put(3,5.5){\circle*{0.1}} 
%\put(3,5.2){\scriptsize{$J$}}
%\put(6,11.2){\scriptsize{$B$}}
%\put(6.2,11.5){\vector(1,2){0.18}}
\put(2.4,3.5){\scriptsize{Binomial}}
\put(1.4,8.8){\scriptsize{Deterministic}}
\put(8.4,8.8){\scriptsize{Deterministic}}
%\put(1.0,8.2){\scriptsize{(area of triangle AGP)}}
%\put(1.7,11.3){\scriptsize{Beta-Pascal (above AP,}}
%\put(1.5,11){\scriptsize{ on the left of halfline PB)}}
%\put(6,7){\scriptsize{Beta-binomial (half-plane below GP)}}
%\put(6.5,11){\scriptsize{Negative binomial (halfline PB)}}
\put(6,3.5){\scriptsize{Deterministic or Binomial or unknown}}
\put(6,3.2){\scriptsize{(depends on the relation between $\alpha$, $\beta$ and $S$)}}

\end{picture}
\end{figure}
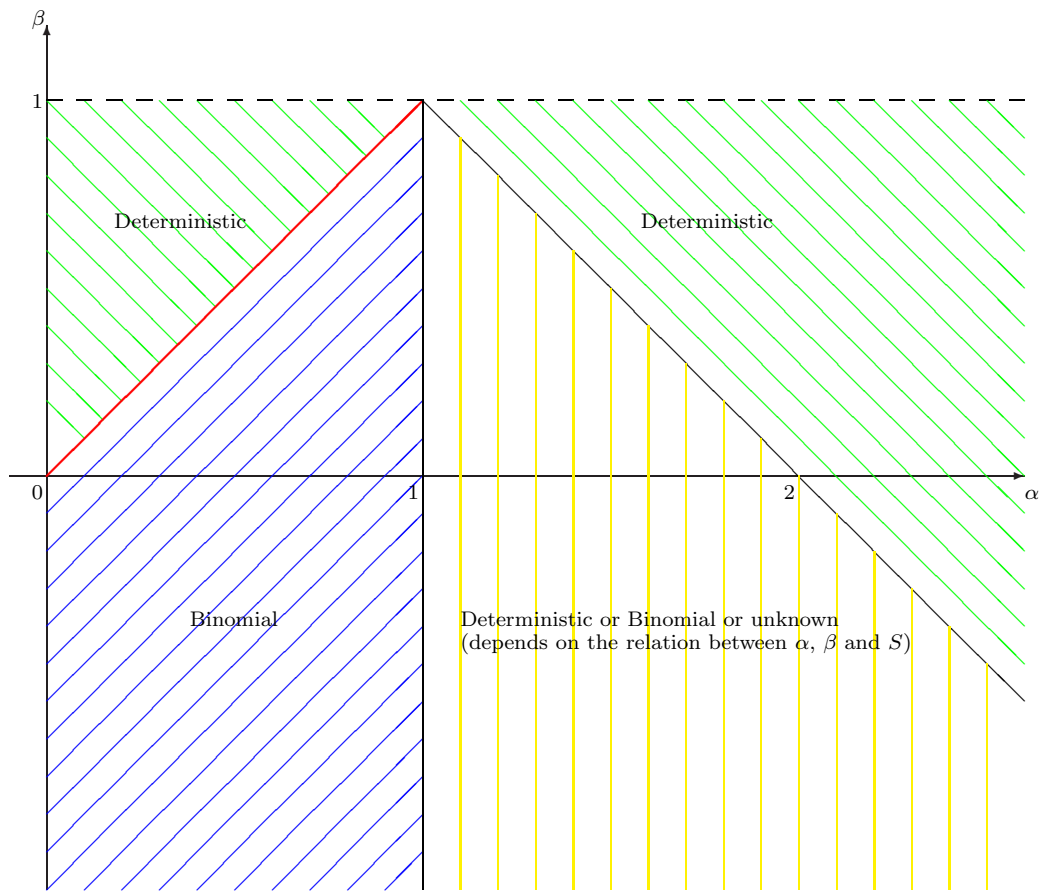

%\section{Concluding remarks}
The power method cannot be applied if $\alpha = \beta$; however, it was shown in \cite{macutek_2006} that in this particular summation the sequence of the descendant distributions converges to the geometric distribution for a wide family of original parents. Specifically, if the original parent is a distribution with a finite support, the limit distribution is deterministic (which can be considered a special case of the geometric distribution).

\section*{Acknowledgement}\label{sec:acknowledgement}
Supported by grants VEGA 2/0054/18 (M. Ko\v{s}\v{c}ov\'a, J. Ma\v{c}utek) and VEGA 1/0341/19 (R. Harman).

%\bibliographystyle{asa}
%\bibliography{clanok}

\end{document}